\documentclass[letterpaper, 10 pt, conference]{ieeeconf}%
\usepackage{subfigure}
\usepackage{bm}
\usepackage{amssymb}
\usepackage{multicol}
\usepackage{multirow}
\usepackage{booktabs} 
\usepackage{epsfig}
\usepackage{amsmath}
\usepackage{graphicx}
\usepackage{epstopdf}
\usepackage{amsfonts}%
\usepackage{indentfirst}
\newtheorem{problem}{\textbf{Problem}}
\newtheorem{definition}{\rm\textbf{Definition}}
\newtheorem{theorem}{\rm\textbf{Theorem}}

\newtheorem{remark}{\rm\textbf{Remark}}
\providecommand{\U}[1]{\protect\rule{.1in}{.1in}}
\IEEEoverridecommandlockouts
\begin{document}
\title{{\LARGE \textbf{Learning Feasibility Constraints for Control Barrier Functions}}}
\author{Wei Xiao, Christos G. Cassandras, and Calin A. Belta \thanks{This work was
		supported in part by NSF under grants ECCS-1931600, DMS-1664644, CNS-1645681, IIS-1723995, and IIS-2024606, by ARPAE under grant DE-AR0001282, by AFOSR under grant FA9550-19-1-0158, and by the MathWorks.  }\thanks{W. Xiao is with the Computer Science and Artificial Intelligence Lab, Massachusetts Institute of Technology \texttt{{\small weixy@mit.edu}}} \thanks{C. G. Cassandras and C. Belta are
with the Division of Systems Engineering and Center for Information and
Systems Engineering, Boston University, Brookline, MA, 02446, USA
\texttt{{\small \{cgc, cbelta\}@bu.edu}}}}
\maketitle
\begin{abstract}
It has been shown that optimizing quadratic costs while stabilizing affine control systems to desired (sets of)
states subject to state and control
constraints can be reduced to a sequence of Quadratic Programs (QPs) by using
Control Barrier Functions (CBFs) and Control Lyapunov Functions (CLFs). In this paper, we employ machine learning techniques to ensure the feasibility of these QPs, which is a challenging problem, especially for high relative degree constraints where High Order CBFs (HOCBFs) are required. 
To this end, we propose a sampling-based learning approach to learn a new feasibility constraint for CBFs; this constraint is then enforced by another HOCBF added to the QPs. The accuracy of the learned feasibility constraint is recursively improved by a recurrent training algorithm. We demonstrate the advantages
of the proposed learning approach to constrained optimal control problems with specific focus on a robot control problem and on autonomous driving in an unknown environment.
\end{abstract}

\section{Introduction}
\label{sec:intro}

With increased interest in autonomous systems, optimal control problems over long finite horizons become increasingly important but challenging, especially in the presence of both safety constraints and control limitations since they may conflict with each other. It was recently shown that for nonlinear control systems that are affine in controls and cost functions that are quadratic in controls, an optimal control problem with safety constraints can be solved through a sequence of quadratic programs (QPs) that are implemented on-line. Central to this approach is the notion of forward invariance enforced using   
barrier functions (BF) \cite{Ames2017}, \cite{Xiao2019}.

BFs are Lyapunov-like functions \cite{Tee2009}, 
\cite{Wieland2007}, whose use can be traced back to optimization problems
\cite{Boyd2004}. More recently, they have been employed to prove set
invariance \cite{Aubin2009}, \cite{Wisniewski2013}. 
Control BFs (CBFs) are extensions of BFs for control systems that are used to
map a constraint defined over system states onto a constraint on the control
input. Recently, it has been shown that, to stabilize an affine control system
while optimizing a quadratic cost and satisfying state and control
constraints, CBFs can be combined with Control Lyapunov Functions (CLFs)
 \cite{Aaron2012} to form
quadratic programs (QPs) \cite{Ames2017}, 
\cite{Glotfelter2017} that are solved in real time.
These CBFs work for constraints
that have relative degree one. A more general form \cite{Nguyen2016}
for arbitrarily high relative degree constraints, termed exponential CBF, employs input-output
linearization and finds a pole placement controller with negative poles. The high order
CBF (HOCBF) proposed in \cite{Xiao2019} is simpler and more general than the
exponential CBF \cite{Nguyen2016}. 

{There are several remaining challenges to be addressed in the CBF-based method, including the determination of (possibly unknown) unsafe sets (e.g., obstacles) and (possibly unknown) system dynamics, as well as the feasibility of the associated QPs when both state constraints (enforced by HOCBFs) and control
bounds are involved. To determine unsafe sets and system dynamics, machine learning techniques are commonly used. Supervised learning techniques have been proposed to learn safe set definitions from demonstrations \cite{Robey2020}, and sensor data \cite{Srinivasan2020}, which are then enforced by CBFs. In \cite{taylor2020learning}, data are used to learn system dynamics for CBFs. In a similar setting, the authors of \cite{Lopez2020} use adaptive control techniques to estimate the unknown system parameters in CBFs.  However, these works did not consider the \emph{feasibility} of the associated QPs in the presence of tight control bounds, which critically depends on how we define CBFs (or HOCBFs), i.e., the parameters involved in the definition of CBFs. In this paper, we focus on the feasibility of these CBF-based QPs.}

Infeasibility can be avoided by precomputing feasible motion spaces \cite{Orthey2013} and using a receding horizon scheme as in Model Predictive Control (MPC) \cite{Diehl2006}.  Optimal control methods can also avoid the infeasibility problem, but are hard to implement with non-linear dynamics and constraints. Furthermore, unknown environments make solving such problems even harder. Some approaches to improve feasibility for specific applications have been proposed. As an example, for the adaptive cruise control (ACC) problem defined in \cite{Ames2017}, the infeasibility issue is addressed by considering the minimum braking distance. In this case, an additional complex safety constraint needs to be added. Further, this approach does not scale well for high-dimensional systems. The penalty method proposed in \cite{Xiao2021TAC2} can improve the recursive feasibility of the QPs and scales well, but it often does not stabilize the system to desired equilibria when ``irregular'' unsafe sets (defined later) are involved. {Feasibility guarantees for CBF-based QPs can also be achieved by finding explicit sufficient conditions \cite{Xiao2021} expressed themselves as CBFs; however, such conditions are usually hard to find for a general constrained control problem.}

The use of machine learning techniques to obtain feasible solutions was recently proposed for legged robots. Feasibility constraints for probabilistic models are learned in \cite{Carpentier2017} based on simplified models. Since the learned constraints are complex, they are simplified by expectation-maximization. Robot footstep limits are modeled as hyperplanes based on success and failure datasets in \cite{Perrin2012}. Reinforcement learning (RL) \cite{Mnih2015} has the potential to address the infeasibility issue for optimal control problems, but it is difficult to quantify infeasibility as a reward; moreover, the optimized parameters may drive the system to a local infeasible region where a feasible solution may never be found.

{The main contribution of this paper is to address the CBF-based QP infeasibility in avoiding all possible unsafe sets in unknown environments using machine learning techniques. Feasibility pertains to the CBF-associated QPs mentioned earlier, which are used in solving a general-purpose optimal control problem. Unsafe sets are assumed to belong to a finite collection of sets whose geometries are known in advance, whereas their locations are only detected during system operation. Our approach proceeds by learning a new feasibility constraint that guarantees feasibility.} 
Specifically, for each type of unsafe set, we sample the state space of the system in its proximity, check for feasibility of the QP one step forward for \emph{regular} unsafe sets and multiple steps forward for \emph{irregular} unsafe sets (precise definitions are provided in the sequel), and learn a differentiable classifier (for feasible and infeasible states) that is then added to the set of initial constraints. We show that, if the initial QP is  feasible, then any control input that satisfies the corresponding CBF constraints renders all the QPs feasible. An additional contribution of the paper is to improve the accuracy of the classifier by introducing a feedback training algorithm that improves the classifier in a recursive way. We validate the effectiveness of the proposed learning-based approach on a robot control problem in an unknown environment that has both regular (e.g., circular) and irregular (e.g., overlapped circular) obstacles, as well as on an application in autonomous driving, in which moving obstacles (such as other vehicles) are usually involved. 


\section{PRELIMINARIES}
\label{sec:pre}

\begin{definition}
	\label{def:classk} (\textit{Class $\mathcal{K}$ function} \cite{Khalil2002}) A
	continuous function $\alpha:[0,a)\rightarrow[0,\infty), a > 0$ is said to
	belong to class $\mathcal{K}$ if it is strictly increasing and $\alpha(0)=0$. A continuous function $\beta:\mathbb{R}\rightarrow\mathbb{R}$ is said to belong to extended class $\mathcal{K}$ if it is strictly increasing and $\beta(0)=0$.
\end{definition}

Consider an affine control system of the form
\begin{equation}
\dot{\bm{x}}=f(\bm x)+g(\bm x)\bm u \label{eqn:affine}%
\end{equation}
where $\bm x\in\mathbb{R}^{n}$, $f:\mathbb{R}^{n}\rightarrow\mathbb{R}^{n}$
and $g:\mathbb{R}^{n}\rightarrow\mathbb{R}^{n\times q}$ are {locally}
Lipschitz, and $\bm u\in U\subset\mathbb{R}^{q}$ with the control constraint set $U$
defined as
\begin{equation}
U:=\{\bm u\in\mathbb{R}^{q}:\bm u_{min}\leq\bm u\leq\bm u_{max}\}.
\label{eqn:control}%
\end{equation}
with $\bm u_{min},\bm u_{max}\in\mathbb{R}^{q}$ and the inequalities are
interpreted componentwise.

\begin{definition}
	\label{def:forwardinv} A set $C\subset\mathbb{R}^{n}$ is forward invariant for
	system (\ref{eqn:affine}) if its solutions {for some $\bm u\in U$} starting at any $\bm x(0) \in C$
	satisfy $\bm x(t)\in C,$ $\forall t\geq0$.
\end{definition}

\begin{definition}
	\label{def:relative} (\textit{Relative degree} \cite{Khalil2002}) The relative degree of a
	(sufficiently many times) differentiable function $b:\mathbb{R}^{n}%
	\rightarrow\mathbb{R}$ with respect to system (\ref{eqn:affine}) is the number
	of times it needs to be differentiated along its dynamics until the control
	$\bm u$ explicitly shows {for all $\bm x$} in the corresponding derivative.
\end{definition}

In this paper, we assume that if there exists $\bm x$ such that the control shows up in the derivative of $b$, then it shows up for all $\bm x$. Since function $b$ is used to define a constraint $b(\bm
x)\geq0$, we will also refer to the relative degree of $b$ as the relative
degree of the constraint. For a constraint $b(\bm x)\geq0$ with relative
degree $m$, $b:\mathbb{R}^{n}\rightarrow\mathbb{R}$, and $\psi_{0}(\bm
x):=b(\bm x)$, we define a sequence of functions $\psi_{i}:\mathbb{R}%
^{n}\rightarrow\mathbb{R},i\in\{1,\dots,m\}$:
\begin{equation}
\begin{aligned} \psi_i(\bm x) := \dot \psi_{i-1}(\bm x) + \alpha_i(\psi_{i-1}(\bm x)),\quad i\in\{1,\dots,m\}, \end{aligned} \label{eqn:functions}%
\end{equation}
where $\alpha_{i}(\cdot),i\in\{1,\dots,m\}$ denotes a $(m-i)^{th}$ order
differentiable class $\mathcal{K}$ function.

We further define a sequence of sets $C_{i}, i\in\{1,\dots,m\}$ associated
with (\ref{eqn:functions}) in the form:
\begin{equation}
\label{eqn:sets}\begin{aligned} C_i := \{\bm x \in \mathbb{R}^n: \psi_{i-1}(\bm x) \geq 0\}, \quad i\in\{1,\dots,m\}. \end{aligned}
\end{equation}

\begin{definition}
	\label{def:hocbf} (\textit{High Order Control Barrier Function (HOCBF)}
	\cite{Xiao2019}) Let $C_{1}, \dots, C_{m}$ be defined by (\ref{eqn:sets}%
	) and $\psi_{1}(\bm x), \dots, \psi_{m}(\bm x)$ be defined by
	(\ref{eqn:functions}). A function $b: \mathbb{R}^{n}\rightarrow\mathbb{R}$ is
	a High Order Control Barrier Function (HOCBF) of relative degree $m$ for
	system (\ref{eqn:affine}) {if there exist} $(m-i)^{th}$ order differentiable
	class $\mathcal{K}$ functions $\alpha_{i},i\in\{1,\dots,m-1\}$ and a class
	$\mathcal{K}$ function $\alpha_{m}$ such that
		\begin{equation}
		\label{eqn:constraint}\begin{aligned} 
		\sup_{\bm u\in U}[L_f^{m}b(\bm x) + [L_gL_f^{m-1}b(\bm x)]\bm u \!+\! O(b(\bm x)) \\+ \alpha_m(\psi_{m-1}(\bm x))] \geq 0, \end{aligned}
		\end{equation}
	for all $\bm x\in C_{1} \cap,\dots, \cap C_{m}$. In
	(\ref{eqn:constraint}), $L_{f}^{m}$ ($L_{g}$) denotes Lie derivatives along
	$f$ ($g$) $m$ (one) times, and $O(b(\bm x)) = \sum_{i = 1}^{m-1}L_f^i(\alpha_{m-i}\circ\psi_{m-i-1})(\bm x).$ {Further, $b(\bm x)$ is such that $L_gL_f^{m-1}b(\bm x)\ne 0$ on the boundary of the set $C_{1} \cap,\dots, \cap C_{m}$.}
\end{definition}

The HOCBF is a general form of the relative degree one CBF \cite{Ames2017},
\cite{Glotfelter2017} (setting $m=1$ reduces the HOCBF to
the common CBF form in \cite{Ames2017}, \cite{Glotfelter2017}), and it is also a general form of the exponential CBF
\cite{Nguyen2016}. We can define $\alpha_i(\cdot)$ in Def. \ref{def:hocbf} to be extended class $\mathcal{K}$ functions to ensure robustness of a HOCBF to perturbations \cite{Ames2017}. 

\begin{theorem}
	\label{thm:hocbf} (\cite{Xiao2019}) Given a HOCBF $b(\bm x)$ from Def.
	\ref{def:hocbf} with the associated sets $C_{1}, \dots, C_{m}$ defined
	by (\ref{eqn:sets}), if $\bm x(0) \in C_{1} \cap,\dots,\cap C_{m}$,
	then any Lipschitz continuous controller $\bm u(t)$ that satisfies
	the constraint in (\ref{eqn:constraint}), $\forall t\geq0$ renders $C_{1}\cap,\dots,
	\cap C_{m}$ forward invariant for system (\ref{eqn:affine}).
\end{theorem}

\begin{definition}
	\label{def:clf} (\textit{Control Lyapunov function (CLF)} \cite{Aaron2012}) A
	continuously differentiable function $V: \mathbb{R}^{n}\rightarrow\mathbb{R}$
	is an exponentially stabilizing control Lyapunov function (CLF) for system
	(\ref{eqn:affine}) if there exist constants $c_{1} >0, c_{2}>0, c_{3}>0$ such
	that for $\forall\bm x\in\mathbb{R}^{n}$, $c_{1}||\bm x||^{2} \leq V(\bm x)
	\leq c_{2} ||\bm x||^{2}, $
	\begin{equation}
	\label{eqn:clf}\underset{u\in U}{inf} \lbrack L_{f}V(\bm x)+L_{g}V(\bm x)
	\bm u + c_{3}V(\bm x)\rbrack\leq0.
	\end{equation}
	
\end{definition}

Many existing works \cite{Ames2017}, \cite{Nguyen2016}, \cite{Yang2019}
combine CBFs for systems with relative degree one with quadratic costs to form
optimization problems. Time is discretized and an optimization problem with
constraints given by the CBFs (inequalities of the form (\ref{eqn:constraint}%
)) is solved at each time step. The inter-sampling effect is considered in \cite{Yang2019}. If convergence to a state is desired, then a
CLF constraint of the form (\ref{eqn:clf}) is added, as in \cite{Ames2017} \cite{Yang2019}. Note that these
constraints are linear in control since the state value is fixed at the
beginning of the interval, therefore, each optimization problem is a quadratic
program (QP) {if the cost is quadratic in the control.} The optimal control obtained by solving each QP is applied at
the current time step and held constant for the whole interval. The state is
updated using dynamics (\ref{eqn:affine}), and the procedure is repeated. Formally, the CBF-based QP is defined as follows. { We 
	partition a time interval $[0,t_f]$  into a set of equal time intervals $\{[0, \Delta t), [\Delta t,2\Delta t),\dots\}$, where $\Delta t > 0$. In each interval $[\omega \Delta t, (\omega+1) \Delta t)$ ($\omega = 0,1,2,\dots$), we assume the control is constant (i.e., the overall control will be piece-wise constant).
	Then {at $t = \omega \Delta t$,} we solve the QP:
	\begin{equation} \label{eqn:obj}
	\begin{aligned}
	\min_{\bm u(\omega \Delta t),\delta(\omega \Delta t)} &\bm u^T(\omega \Delta t) H \bm u(\omega \Delta t) + p_0\delta^2(\omega \Delta t)\\
	&\text{s.t. }\bm u_{min}\leq\bm u\leq\bm u_{max}\\
	&L_fV(\bm x)+L_gV(\bm x) \bm u + \epsilon V(\bm x) \leq \delta,\\
	L_f^{m}b(\bm x) + [L_g&L_f^{m-1}b(\bm x)]\bm u \!+\! O(b(\bm x)) + \alpha_m(\psi_{m-1}(\bm x)) \geq 0
	\end{aligned}
	\end{equation}
	 In the above equation, $H$ is positive definite, $\delta(t)$ is a relaxation variable on the CLF constraint used to avoid conflict with the HOCBF constraint, and $p_0 > 0$ is a penalty on the relaxation $\delta(t) \in\mathbb{R}$. } This method works conditioned on the
fact that the QP at every time step is feasible. However, this is not
guaranteed, in particular under tight control bounds (or very limited controls). In this paper, we show
how the QP feasibility can be recursively improved by using machine learning techniques.

\section{Problem Formulation and Approach}
\label{sec:prob}

Consider an optimal control problem for system (\ref{eqn:affine}) with the cost defined as:
\begin{equation}\label{eqn:gcost}
\int_{0}^{t_f}\mathcal{C}(||\bm u(t)||) dt + p_0||\bm x(t_f) - \bm K||^2,
\end{equation}
where $||\cdot||$ denotes the 2-norm of a vector; $t_f$ denotes a given final time; and $\mathcal{C}$ is a strictly increasing function of its argument (usually quadratic). $\bm K\in\mathbb{R}^n$ is an equilibrium for system (\ref{eqn:affine}) {in the absence of control} and $p_0 > 0$.

$\textbf{Constraint 1}$ (Unsafe state sets): Let $S$ denote an index set for unsafe (state) sets. 
System (\ref{eqn:affine}) avoids each unsafe set $j\in S$ if the state of system (\ref{eqn:affine}) satisfies:
\begin{equation} \label{eqn:obstacle}
b_j(\bm x(t))\geq 0, \forall t\in[0,t_f], {\forall j\in S},
\end{equation}
where $b_j: \mathbb{R}^n\rightarrow\mathbb{R}$ is a continuously differentiable function (not a CBF or HOCBF yet).


$\textbf{Constraint 2}$ (State and control limitations): Assume we have a set of constraints on control input of system (\ref{eqn:affine}) as in (\ref{eqn:control}) and on the state in the form:
\begin{equation}\label{eqn:state}
\begin{aligned}
\bm x_{min}\leq \bm x(t)\leq \bm x_{max}, \forall t\in[0,t_f]
\end{aligned}
\end{equation}
where $\bm x_{min}\in\mathbb{R}^n$ and $\bm x_{max}\in\mathbb{R}^n$ denote the minimum and maximum
state vectors respectively, and the inequalities are interpreted componentwise. Note that system state constraints can usually be relaxed (and the relaxation is minimized, as shown in the CLF in (\ref{eqn:obj})), and therefore are distinguished from the safety constraint (\ref{eqn:obstacle}). The control bounds in (\ref{eqn:control}) usually denote the system control capability, and thus are hard constraints.

A control policy for system (\ref{eqn:affine}) is $\bm {feasible}$ if the hard constraints (\ref{eqn:obstacle}) and (\ref{eqn:control}) are satisfied. 

In this paper, we consider the following problem:

\vspace{1mm}
\begin{problem}\label{prob:main}
	Find a $\bm {feasible}$ control policy for system (\ref{eqn:affine}) such that cost (\ref{eqn:gcost}) is minimized, and state constraints (\ref{eqn:state}) are satisfied.
\end{problem}

{\textbf{Approach:} Problem \ref{prob:main} is a general problem definition and is assumed to be feasible by itself. In order to efficiently solve it, we replace the safety constraints above by HOCBF constraints. This method is conservative as the satisfaction of the CBF constraint is only a sufficient condition for the satisfaction of the original constraint; hence, the resulting solution is sub-optimal but provides safety guarantees. This solution can be driven to near-optimality by optimally tracking a desired trajectory through trajectory planning methods \cite{Xiao2021B}.
} 
In our approach, we use HOCBFs \cite{Xiao2019} to enforce the safety constraint (\ref{eqn:obstacle}) and a CLF \cite{Aaron2012} to enforce the terminal cost in (\ref{eqn:gcost}).

\begin{figure}[thpb]
	\centering
	\includegraphics[scale=0.28]{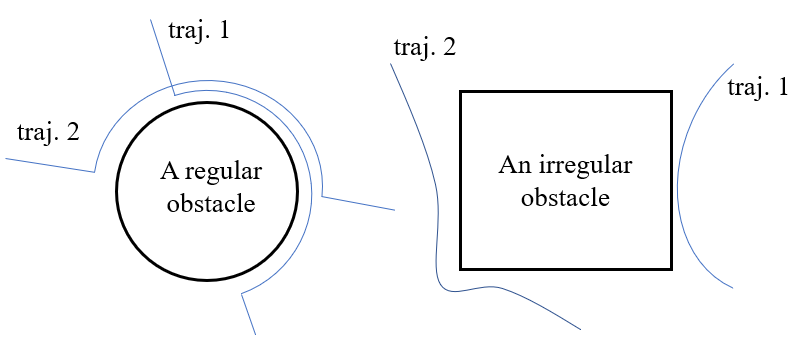}
	\vspace{-3mm}
	\caption{Regular and irregular obstacle (unsafe set) examples for a robot with nonholonomic dynamics. The robot needs the same control input effort (from the CBF-based QPs) in order to avoid the regular circular obstacle, regardless of where the robot is initially located, as shown in example trajectories 1 and 2. However, the robot needs larger control input effort for trajectory 2 than the one for trajectory 1, as there are corners in the irregular rectangular obstacle. Therefore, the feasibility of a CBF-based QP is indeed dependent on the robot initial condition (location).}
	\label{fig:unsafe_sets}%
	\vspace{-3mm}
\end{figure}

The approach to Problem \ref{prob:main} proposed in \cite{Ames2017} is based on 
partitioning the time interval $[0, t_f]$ into $[t_k, t_{k+1}), k = \{0, 1, 2,\dots\}, t_0 = 0$, as introduced at the end of Sec. \ref{sec:pre}.
Since the state is kept constant at its value at $t_k$, {a HOCBF constraint is linear in control, thus, the optimization problem is a QP if the cost is quadratic in the control at $t_k$}. Such a QP can easily become infeasible since (\ref{eqn:control}) may conflict with the HOCBF constraints corresponding to (\ref{eqn:obstacle}). {Depending on how the system initial state may affect the feasibility of the CBF-based QPs, we classify unsafe sets into two classes:
\begin{definition} \label{def:regular}
	(Regular and irregular unsafe sets) Assume Problem \ref{prob:main} is feasibile. An unsafe set $C_u:=\{\bm x\in\mathbb{R}: b(\bm x) < 0\}$ considered in the QPs (\ref{eqn:obj}) is defined as {\it regular} if the feasibility of all the CBF-based QPs (\ref{eqn:obj}) does not depend on the initial state $\bm x(0)$ of system (\ref{eqn:affine}). Otherwise, we say that the set is {\it irregular}. 
\end{definition}

{When there are multiple unsafe sets, whether each one is regular or not is checked one by one through the QP (\ref{eqn:obj}); these can then be combined into a single QP.} An irregular unsafe set generally depends on the dynamics (\ref{eqn:affine}) and corresponds to irregular shapes, such as unsafe sets with sharp corners, in which case the system requires (locally) large control input from the CBF-based QP to avoid corners if the dynamics are nonholonomic; however, the CBF-based QPs may still be feasible if the system trajectory never approaches a corner. An example of a regular unsafe set is a circular obstacle, and an example of an irregular unsafe set is a rectangle for a robot with nonholonomic dynamics, as shown in Fig. \ref{fig:unsafe_sets}.}

Moreover, In order to address the infeasibility problem of the CBF-based QPs in an offline way, 
we first define unsafe sets as being of the same ``type'' if they have the same geometry, meaning the conditions for problem feasibility are the same, e.g., circular unsafe sets are the same type if they have the same radius but different locations. {Let $S_t$ denote the  set indexing all the unsafe set types.} In an unknown environment, if unsafe set types are known, then we can study the problem feasibility based on these types. Otherwise, we can study the problem feasibility based on some selected set types, and then use these types of unsafe sets to over-approximate other unknown types; for example, a regular circular unsafe set can be used to over-approximate any shape of unsafe sets. In this paper, we limit ourselves to a set of known unsafe set types (a typical application is in autonomous driving where the vehicle types are known). 

We propose a learning-based approach, specifically a classifier, to ensure the QP feasibility for a certain type of unsafe set. For each type of unsafe set, we sample states in the vicinity of the set. For each of the sampled states, we solve the QP (one or more steps forward, and more steps can help to check the feasibility of a receding horizon) and label the state as +1 if the QP is feasible for all steps, and label it as -1 otherwise. {In contrast to parameter learning,  sampling-based learning is more myopic as the parameter learning approach considers the horizon from the initial time to the final time.} Then, we use a machine learning algorithm to classify all the feasible and infeasible states, and get a feasibility constraint from the classifier hypersurface. This feasibility constraint is then enforced by a HOCBF and added to the QP. This learning approach allows us to deal with irregular unsafe sets in an unknown environment. 

\section{Sampling-Based Learning Approach}
\label{sec:learn}

In this section, we show how we can deal with irregular unsafe sets as defined in Def. \ref{def:regular}. This approach also works for regular unsafe sets, but tends to be conservative. Recall that the type of every unsafe set $i\in S$ is already known in an unknown environment, $S_t$ denotes an index set for unsafe set types in an unknown environment, and $S_j\subseteq S, j\in S_t$ denotes the index set for unsafe sets of type $j$.


\subsection{Feasible and Infeasible State Sets}
\label{sec:FI}

The QP (\ref{eqn:obj}) may be infeasible at a given state $\bm x(t)$ at time $t$. The constraints in (\ref{eqn:obstacle}) form a constraint set for the state of system (\ref{eqn:affine}). Without control (i.e., $\bm u(t) = 0, \forall t\in[0,t_f]$), system (\ref{eqn:affine}) may escape from this constraint set for a given initial state $\bm x(0)$. However, if the system is controlled with the optimal control $\bm u^*(t)$ from solving the QP (\ref{eqn:obj}), the system may also exit this constraint set since the limited control may not be able to prevent the system from leaving this set when the state approaches the set boundary, which typically happens in high relative degree systems. Then, QP (\ref{eqn:obj}) becomes infeasible.

An intuitive example is a robot control problem. Suppose a robot, with limited  control input (deceleration), needs to arrive at a destination while avoiding an obstacle that is located between the robot's initial position and the destination. When the robot gets close to the obstacle with high speed, it may not be able to brake in time to avoid the obstacle since the control in QP (\ref{eqn:obj}) is limited by (\ref{eqn:control}). However, when the speed is low, the robot can safely avoid the obstacle. 


The main idea of the sampling-based learning approach is to partition the state space of system (\ref{eqn:affine}) into sets in which the QP (\ref{eqn:obj}) is feasible or infeasible after a certain number of time steps. 
This is a difficult problem, especially for high-dimensional systems with fast dynamics. We show how machine learning techniques can be used to address this.

\subsection{Sampling and Classification} 
\label{sec:classify}
As mentioned before, the system is in an unknown environment such that it only knows the types of the unsafe sets the environment may include, but not their number and locations.
In order to make the learned feasibility constraint independent from the location of an unsafe set, we choose the relative coordinate $\bm z\in\mathbb{R}^n$ between the system and unsafe set as one of the input features for machine learning techniques. For example, let the system state be $\bm x := (x_1, x_2,\dots, x_n)$. If $x_1, x_2$ denote the 2-D position of an object in $\bm x$, then we define input features $\bm z:= (x_1 - x_o, x_2 - y_o, x_3, \dots, x_n)$ for the machine learning techniques, where $(x_o,y_o)\in \mathbb{R}^2$ denotes the 2-D location center
of the unsafe set. Along the same lines, we may also consider the relative speed and acceleration between the system and unsafe set as the input for the machine learning model in order to consider moving unsafe sets.

For each type of unsafe set $j\in S_t$, since we only consider the relative coordinates as the input for the learning model as discussed above, we arbitrarily assign a location and an orientation (if it exists) for $j$ and randomly sample around $j$ to find an initial state $\bm z(0)$ around the unsafe set. We then solve the QP (\ref{eqn:obj}) at time $0$ according to the geometry of the unsafe set:
\begin{itemize}
	\item (i) \textbf{Regular unsafe set:} we solve the QP (\ref{eqn:obj}) at time $0$ for $1$ time step forward.
	\item (ii) \textbf{Irregular unsafe set:} we solve the QP (\ref{eqn:obj}) at time $0$ for $H_t \in \mathbb{N} > 1$ time steps forward.
\end{itemize}

\begin{remark}
	Unlike regular unsafe sets, when dealing with irregular unsafe sets, the system may get stuck at local traps. This is why we extend the solution of the QP (\ref{eqn:obj}) to $H_t > 1$ time steps. In this case, any one of the $H_t$-step QPs becoming infeasible will make the system fail. The local traps can easily make the QP (\ref{eqn:obj}) infeasible, especially when the system state approaches their boundary. Therefore, it is more likely to make an inital state that is located around the local traps belong to the infeasible set when we solve the QP (\ref{eqn:obj}) $H_t > 1$ steps forward. Then, the system may avoid the local traps if it avoids the infeasible set, and thus improve its reachability.
\end{remark}

If the QP (or all the QPs in case (ii)) (\ref{eqn:obj}) is feasible, we label the state $\bm z(0)$ as $+1$. Otherwise, it is labelled as $-1$. This procedure results in two labelled classes. We employ a machine learning technique (such as Support Vector Machine (SVM), Deep Neural Network (DNN), etc.) with $\bm z(0)$ as input to perform classification, and get a classification hypersurface for each $j\in S_t$ in the form:
\begin{equation}\label{eqn:hyperS}
H_j(\bm z):\mathbb{R}^n\rightarrow\mathbb{R},
\end{equation}
where $H_j(\bm z(0))\geq 0$ denotes that $\bm z(0)$ belongs to the feasible set. This inequality is called the \textbf{feasibility constraint}. 

 Assuming the relative degree of (\ref{eqn:hyperS}) is $\gamma \in \mathbb{N}$, we define the set of all control values that satisfy $H_j(\bm z(t))\geq 0$ as:
\begin{equation}\label{eqn:feaset}
\begin{aligned}
K_{fea}^j = \{\bm u\in U: L_f^{\gamma}H_j(\bm z) + L_gL_f^{\gamma-1}H_j(\bm z)\bm u  \\+ O(H_j(\bm z))  + \alpha_{\gamma}(\psi_{\gamma-1}(\bm z)) \geq 0\}
\end{aligned}
\end{equation}
where $\psi_{m-1}$ is recursively defined as in (\ref{eqn:functions}) by $H_j$ with extended class $\mathcal{K}$ functions. 

We define feasibility forward invariance as follows:

\begin{definition} \label{def:feainv}
	An optimal control problem is feasibility forward invariant for system (\ref{eqn:affine}) if its solutions starting at all feasible $\bm x(0)$  are feasible for all $t\geq 0$.
\end{definition} 

\begin{theorem} \label{thm:fea}
	Assume that the hypersurfaces $H_j(\bm z), \forall j\in S_t$ ensure $100\%$ feasibility and infeasibility classification accuracy. If $H_j(\bm z(0))\geq 0, \forall j\in S_t$, then any Lipschitz continuous controller $\bm u(t)\in K_{fea}^j, \forall j\in S_t$ renders Problem \ref{prob:main} feasibility forward invariant.
\end{theorem}

\textbf{Proof:} By Theorem 5 in \cite{Xiao2019} and $H_j(\bm z(0))\geq 0, \forall j\in S_t$, any control input that satisfies $\bm u(t)\in K_{fea}^j, \forall j\in S_t, \forall t\in[0, \infty]$ makes $H_j(\bm z(t))\geq 0, \forall j\in S_t, \forall t\in[0, \infty]$. Since $H_j(\bm z), \forall j\in S_t$ classifies the state space of system (\ref{eqn:affine}) into feasible and infeasible spaces for Problem \ref{prob:main} with $100\%$ accuracy, we have that Problem \ref{prob:main} is feasibility forward invariant for system (\ref{eqn:affine}).$\qquad\qquad\qquad\qquad\qquad\qquad\qquad\blacksquare$

Naturally, machine learning techniques cannot ensure $100\%$ classification accuracy. We introduce an approach 
based on feedback training to improve the accuracy in the following subsection. In fact, if the classification accuracy is high enough, Problem \ref{prob:main} may also be always feasible since system (\ref{eqn:affine}) may never reach the infeasible space.

Similar to the QP (\ref{eqn:obj}), we have a feasible reformulated problem at $t = \omega \Delta t$ ($\omega = 0,1,2,\dots, \frac{t_f}{\Delta t}-1$):
{\small\begin{equation} \label{eqn:objf}
\begin{aligned}
\min_{\bm u(t),\delta(t)}& \bm u^T(t) H \bm u(t) + p_0\delta^2(t)\\
&\text{s.t. }\bm u_{min}\leq\bm u\leq\bm u_{max}\\
&L_fV(\bm x)+L_gV(\bm x) \bm u + \epsilon V(\bm x) \leq \delta,\\
L_f^{m}b(\bm x) + &[L_gL_f^{m-1}b(\bm x)]\bm u \!+\! O(b(\bm x)) + \alpha_m(\psi_{m-1}(\bm x)) \geq 0,\\
L_f^{\gamma}H_j(\bm z) +& L_gL_f^{\gamma-1}H_j(\bm z)\bm u  + O(H_j(\bm z)) + \alpha_{\gamma}(\psi_{\gamma-1}(\bm z)) \geq 0
\end{aligned}
\end{equation}
}where $b(\bm x) = b_j(\bm x)$, and every safety constraint of the same type $j$ uses the same $H_j(\bm z(t))\geq 0, \forall j\in S_t$. 

\subsection{Feedback Training} 
\label{sec:train}

For each $j\in S_t$, we first sample the points without any hypersurface (\ref{eqn:hyperS}). After the first iteration, we obtain a hypersurface that classifies the state space of system (\ref{eqn:affine}) into feasible and infeasible sets, but with relatively low accuracy. Then we can add these hypersurfaces (\ref{eqn:hyperS}) into the QP (\ref{eqn:obj}) (i.e., use (\ref{eqn:objf})) and sample new data points to perform a new classification, and obtain another classification hypersurface that replaces the old one. Iteratively, the classification accuracy is improved and the infeasible set shrinks. 

\begin{figure}[thpb]
	\centering
	\vspace{-3mm}
	\includegraphics[scale=0.25]{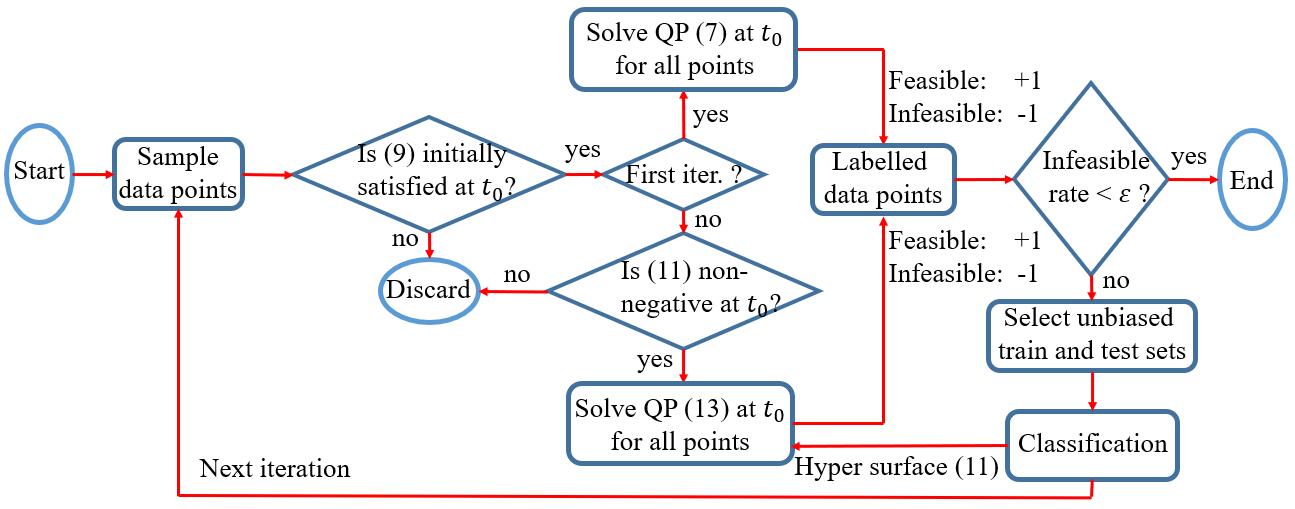}
	\caption{Feedback training workflow for unsafe set $j\in S_t$ ($\varepsilon > 0$ denotes the termination threshold).}
	\label{fig:work}%
	\vspace{-3mm}
\end{figure}

Since the CBF method requires the constraint to be initially satisfied, we discard the samples that do not meet this requirement. To ensure classification accuracy, we also need unbiased data samples. The workflow is shown in Fig. \ref{fig:work}. The infeasibility rate is the ratio of the number of infeasible samples over the total times of solving the QP (\ref{eqn:obj}) or (\ref{eqn:objf}).

\subsection{Generalization}
\label{sec:tech}

Since we sample data around the unsafe set, we also need to check the generalization of the hypersurface (\ref{eqn:hyperS}) in the area where we do not sample since system (\ref{eqn:affine}) may actually start from some state in the unsampled area. Problem \ref{prob:main} is usually feasible when system (\ref{eqn:affine}) is far away from the unsafe set, therefore, the unsampled area should be located at the positive side of the hypersurface (\ref{eqn:hyperS}), which can be viewed as the generalization of this hypersurface, as usually appears in machine learning techniques. 

Once we get a hypersurface (feasibility constraint) for a type of unsafe set in the pre-training process, we can also apply this feasibility constraint to other unsafe sets that are of the same type but with different locations since the hypersurface only depends on the relative location of the pre-training unsafe set. This is helpful for systems in which we do not know the number and locations of the unsafe set, but know the type of unsafe sets the environment has.

It is important to note that the optimal hypersurface is not unique given the training samples; this is due to the weight space symmetries \cite{Bishop2006} in neural networks.

\section{Implementation and Case Studies}
\label{sec:sim}

We implemented the proposed learning approach in MATLAB and performed simulations for a robot control problem. Suppose all the obstacles are of the same type but the obstacle number and their locations are unknown to the robot, and the robot is equipped with a sensor ($\frac{2}{3}\pi$ field of view (FOV) and $7m$ sensing distance with $1m$ uncertainty) to detect the obstacles.

With $\bm x := (x,y,\theta,v), \bm u = (u_1,u_2)$, the dynamics are defined as:
\begin{equation}\label{eqn:robot}
\begin{aligned}
\dot x = v\cos\theta,\; \;\;\dot y = v\sin\theta, \; \;\;
\dot \theta = u_1, \;\;\dot v = u_2,
\end{aligned}
\end{equation}
where $x, y$ denote the location along $x, y$ axis, respectively, $\theta$ denotes the heading angle of the robot, $v$ denotes the linear speed, and $u_1, u_2$ denote the two control inputs for turning and acceleration, respectively. 

{We instantiate cost (\ref{eqn:gcost}) in the form:
\begin{equation}
\min_{\bm u(t)} \int_{0}^{t_f} \left[u_1^2(t) + u_2^2(t)\right] dt + p_0((x(t_f) - x_d)^2 + (y(t_f) - y_d)^2).
\end{equation}
In other words, we wish to minimize the energy consumption and drive the robot to a given destination $(x_d, y_d)\in\mathbb{R}^2$,} i.e., drive $(x(t),y(t))$ to $(x_d, y_d), \forall t\in[t',t_f], \text{for some }t'\in[0,t_f]$, as defined in (\ref{eqn:gcost}). The robot dynamics are not full state linearizable \cite{Khalil2002} and the relative degree of the position (output) is 2. Therefore, we cannot directly apply a CLF. However, the robot can arrive at the destination if its heading angle $\theta$ stabilizes to the desired direction and its speed $v$ stabilizes to a desired speed $v_0 > 0$, i.e.,
\begin{equation}\label{eqn:robotdst}
\theta(t) \rightarrow \arctan(\frac{y_d - y(t)}{x_d - x(t)}),\;\; 
v(t)\rightarrow v_0, \forall t\in[0,t_f].
\end{equation}
Now, we can apply the CLF method since the relative degrees of the heading angle and speed are 1.

The unsafe sets (\ref{eqn:obstacle}) are defined as circular (regular) obstacles:
\begin{equation}\label{eqn:robotobs}
\sqrt{(x(t) - x_i)^2 + (y(t) - y_i)^2} \geq r, \forall i\in S,
\end{equation}
where $(x_i, y_i)$ denotes the location of the obstacle $i\in S$, and $r > 0$ denotes the safe distance to the obstacle. Note that we may have irregular obstacles when there are overlapped circular obstacles.

The speed and control constraints (\ref{eqn:control}) are defined as:
\begin{equation}
\begin{aligned}
V_{min} \leq v(t) \leq V_{max},
\end{aligned}
\end{equation}
\begin{equation} \label{eqn:robotcontrol1}
u_{1,min} \leq u_1(t) \leq u_{1,max}, u_{2,min} \leq u_2(t) \leq u_{2,max},
\end{equation}
where $V_{min} = 0m/s, V_{max} = 2m/s, u_{1, max} = -u_{1, min} = 0.2rad/s, u_{2, max} = -u_{2, min} = 0.5m/s^2$. Other parameters are $p_0 = 1, \Delta t = 0.1s, \epsilon = 10$.  

\subsection{Robot Control}

We chose all class $\mathcal{K}$ functions in the definitions of all HOCBFs as linear functions, and used SVM to classify the feasible and infeasible sets for the QP (\ref{eqn:obj}) or (\ref{eqn:objf}). We obtained a hypersurface for each type of obstacle, and apply this hypersurface to the same type obstacles (both regular and irregular) with unknown locations.

\subsubsection{Regular Obstacles}

\begin{figure}[htbp]
	\vspace{-3mm}
	\centering
	\subfigure[Trajectories]{
		\begin{minipage}[t]{0.21\textwidth}
			\centering
			\includegraphics[width=\textwidth]{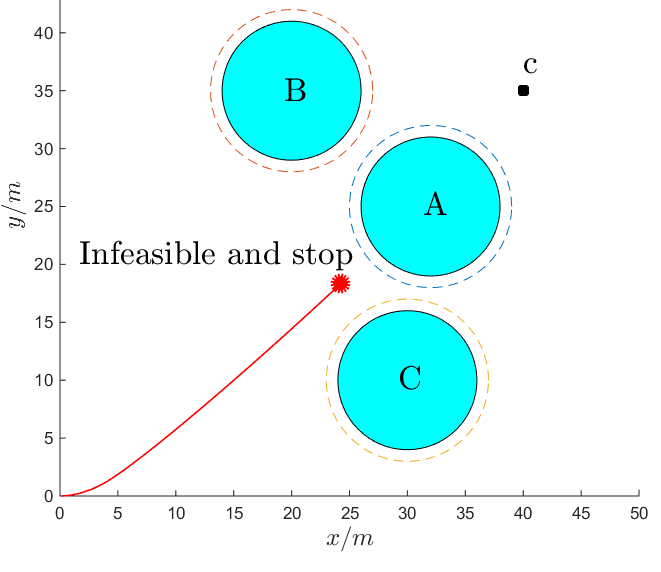}
		\end{minipage}\label{fig:infa}%
	}
	\subfigure[Control $u_1(t)$]{
		\begin{minipage}[t]{0.24\textwidth}
			\centering
			\includegraphics[width=\textwidth]{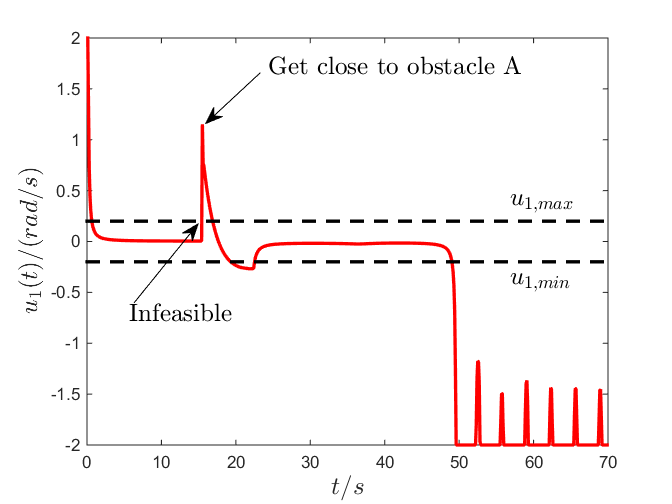}
		\end{minipage}\label{fig:infb}%
	}\hfill	

	
	\centering
	\caption{The control profiles for the infeasible example with relaxation on both control limitations (\ref{eqn:robotcontrol1}). The control bound for $u_2(t)$ is satisfied.}\label{fig:infea}
\end{figure}

We consider three regular obstacles (locations unknown to the robot, but detected by the on-board sensor with sensing range $> 6m$ in radius), and solve the QP (\ref{eqn:obj}) or (\ref{eqn:objf}) one time step forward to check the feasibility of the samples. Some other parameters are: $r = 7m, \varepsilon = 0.001, \xi = 1m, x_A = 32m, y_A = 25m, x_B = 20m, y_B = 35m, x_C = 30m, y_C = 10m$, and the map of the environment is shown in Fig. \ref{fig:infa}.

Assume the robot destination is $(40m,35m)$. If the QP (\ref{eqn:obj}) is solved with all the objectives and constraints defined as (\ref{eqn:robot})-(\ref{eqn:robotcontrol1}), it will become infeasible when the robot gets close to the obstacle, as shown in Fig. \ref{fig:infa}. In order to show how these two control inputs may lead to the infeasibility of the QP (\ref{eqn:obj}), we relax both limitations (\ref{eqn:robotcontrol1}) and show the control profile in Fig. \ref{fig:infb}.

To solve this infeasiblity problem, we apply the learning method introduced in Sec. \ref{sec:learn}. During pre-training, we arbitrarily assign the location $(x_o,y_o):=(20m,35m)$ of an obstacle that is the same type (circle with the same radius) as $j\in S_t$. Then we define $\bm z:=(x - x_o, y - y_o,\theta, v)$ as input for SVM with polynomial kernel of degree 2, i.e., the kernel $k(\bm y,\bm z)$ is defined as:
\begin{equation}\label{eqn:kernel}
k(\bm y,\bm z) = (k_1 + k_2\bm y ^T \bm z)^2.
\end{equation}
where $\bm y$ denotes an input vector similar to $\bm z$, $k_1\in\mathbb{R}, k_2\in\mathbb{R}$.

\begin{figure}[thpb]
	\centering
	\vspace{-3mm}
	\includegraphics[scale=0.35]{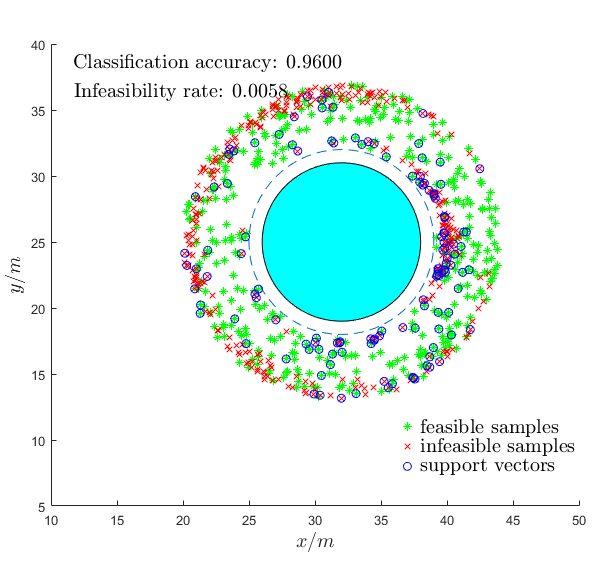}
	\caption{Feedback training at the 3rd iteration for the regular obstacle. All data are sampled around the obstacle (solid circle in all sub-figures). Each sample is a four dimensional point, but is visualized in $x-y$ plane.}
	\label{fig:trainS}%
	\vspace{-3mm}
\end{figure}

We set $k_1 = 0.9, k_2 = 0.4$ for the kernel (\ref{eqn:kernel}) in the feedback training process. The 3rd training iteration for the regular obstacle is shown in Fig. \ref{fig:trainS}. 

\begin{table}
	\caption{Training results for the regular obstacle}
	\label{tab:accu}
	\centering
	\begin{tabular}{l|l|l|ll}
		\toprule
		iter.     & \multicolumn{1}{|c|}{QP (\ref{eqn:objf}) inf. rate}  &  \multicolumn{1}{|c|}{classi. accu.}   & \multicolumn{1}{|c}{train}     & \multicolumn{1}{c}{test}  \\
		\midrule
		1  & 0.0665  &0.8490  & \multicolumn{1}{c}{2k}    & \multicolumn{1}{c}{1k(20k)}   \\
		2  &0.0528& 0.9150   &\multicolumn{1}{c}{0.8k}  & \multicolumn{1}{c}{0.2k (10k)}  \\
		3 & 0.0058& 0.9600&\multicolumn{1}{c}{0.6k} & \multicolumn{1}{c}{0.2k (60k)} \\
		test$^1$ &0.0004&  &   &\multicolumn{1}{c}{25k}\\
		gen.$^2$ &0  &1.0000 &  &  \multicolumn{1}{c}{20k}\\
		\bottomrule
	\end{tabular}
	
	$^1$ denotes testing results within the training data sampling space.
	
	$^2$ denotes testing results out of the training data sampling space.
	\vspace{-5mm}
\end{table}

\begin{figure}[htbp]
	\centering
	\vspace{0mm}
	\subfigure[Trajectories]{
		\begin{minipage}[t]{0.135\textwidth}
			\centering
			\includegraphics[width=\textwidth]{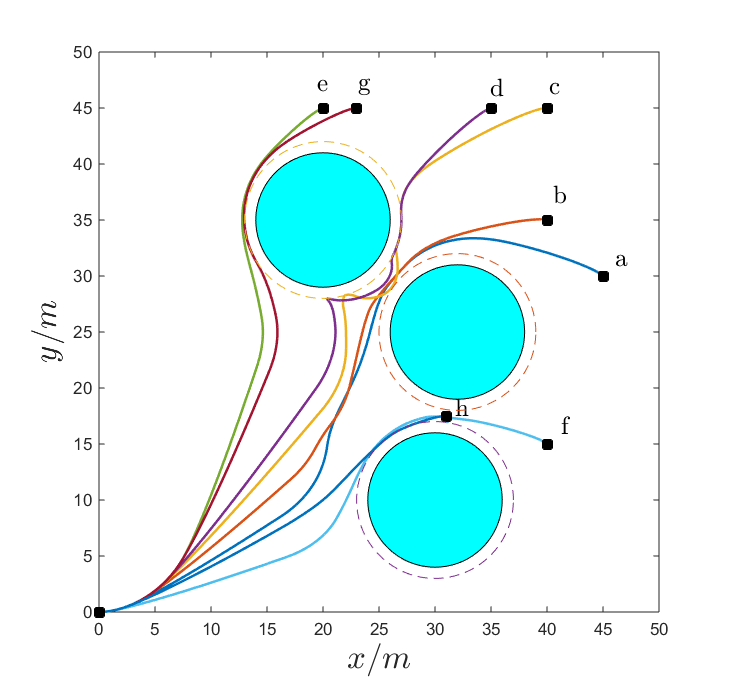}
		\end{minipage}%
	}%
	\subfigure[Control $u_1(t)$]{
		\begin{minipage}[t]{0.17\textwidth}
			\centering
			\includegraphics[width=\textwidth]{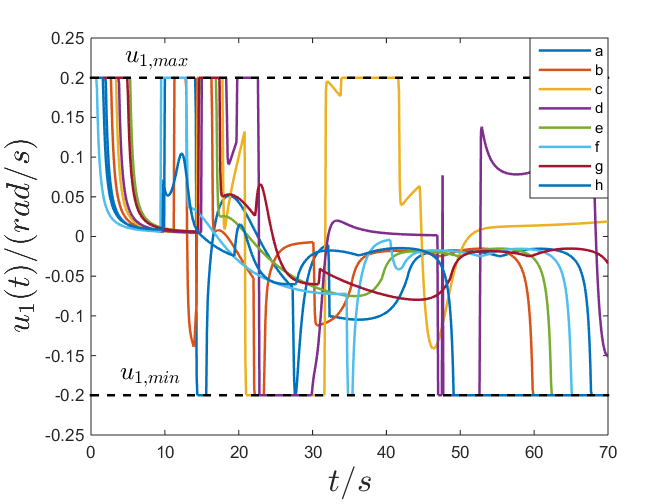}
		\end{minipage}%
	}%
	\subfigure[Control $u_2(t)$]{
		\begin{minipage}[t]{0.17\textwidth}
			\centering
			\includegraphics[width=\textwidth]{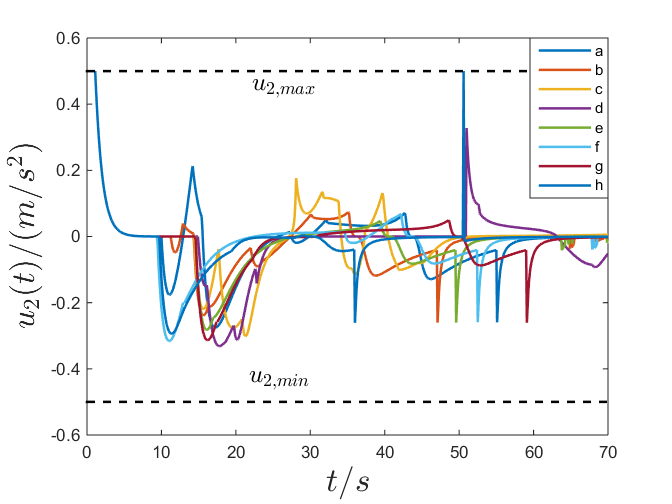}
		\end{minipage}
	}%
	\centering
	\caption{Robot control problem feasibility test after learning in an environment with regular obstacles.}\label{fig:fea}
\end{figure}

As shown in Table \ref{tab:accu}, the classification accuracy for the infeasible and feasible samples is improved after each iteration, and the infeasibility rate also decreases. It becomes inefficient to get infeasible samples after just three iterations since we want to get an unbiased training data set and the infeasibility rate is $0.0004$ (improved 166 times compared to the first iteration). The training results are shown in Table \ref{tab:accu} (the first number of test samples denotes the test number for the classification accuracy, and the number in the bracket denotes the test number for the infeasibility rate). Eventually, we get a hypersurface $H(\bm z)$ from the kernel.

\textbf{Hypersurface generalization:} A natural question is what happens if the robot starts at points outside the sample area (we only sampled around the obstacle with location radius within $[7,13]$, and achieved 96.00\% classification accuracy over the test samples and 0.04\% infeasibility rate for the QP (\ref{eqn:objf})). Therefore, we tested the classification accuracy for samples with location radius within $[13,32]$ (out of the training and test sets). The test accuracy is  \textbf{{ $\bm{100\%}$} (i.e., $H(\bm z) \geq 0$) and the infeasibility rate of the QP (\ref{eqn:objf}) is { $\bm{0\%}$}} (over 20000 samples) for the obstacle (listed in the last row of Table \ref{tab:accu}), which shows good generalization.

We apply this hypersurface to the robot control problem with the objective and constraints defined as (\ref{eqn:robot})-(\ref{eqn:robotcontrol1}) and test 8 different destinations a, b, c, d, e, f, g, h shown in Fig. \ref{fig:fea}. The initial conditions at time $0s$ are $\bm x(0) = (0m,0m,0 rad,1m/s)$ (out of the training and test sets), and $ t_f = 70s$. The QPs (\ref{eqn:objf}) are always feasible during $[0s, 70s]$. The obstacles are safely avoided, and the robot eventually arrives at the destination. We present the two control input profiles in Fig. \ref{fig:fea}.

\subsubsection{Irregular Obstacle}

\begin{figure}[thpb]
	\centering
	\vspace{-3mm}
	\includegraphics[scale=0.35]{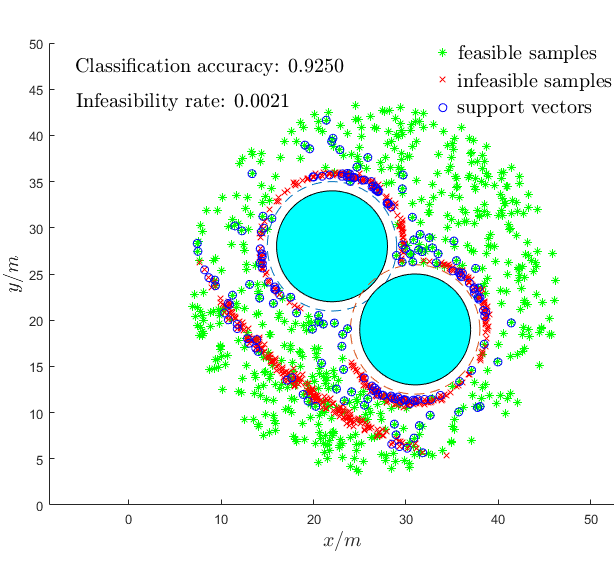}
	\caption{Feedback training at the 3rd iteration for the irregular obstacle. All data are sampled around the obstacle (solid circle in all sub-figures). Each sample is a four dimensional point, but is visualized in $x-y$ plane.}
	\label{fig:train2}%
	\vspace{-3mm}
\end{figure}

In this case, we consider an irregular obstacle that is formed by two  overlapped disks (with locations $(22m,28m)$ and $(31m, 19m)$ but unknown to the robot), as shown in Fig. \ref{fig:train2}. We apply the learning method introduced in Sec. \ref{sec:learn} to recursively improve the problem feasibility, and possibly to escape from local traps. We formulate a receding horizon control of $H_t = 60$, and check the feasibility of all these $H_t$ step QPs. The other settings are the same as the ones from the regular case. The 3rd training iteration is shown in Fig. \ref{fig:train2}.

\begin{table}
	\caption{Training results for the irregular obstacle}
	\label{tab:accu2}
	\centering
	\begin{tabular}{l|l|l|ll}
		\toprule
		iter.     & \multicolumn{1}{|c|}{QP (\ref{eqn:objf}) inf. rate}  &  \multicolumn{1}{|c|}{classi. accu.}   & \multicolumn{1}{|c}{train }     & \multicolumn{1}{c}{test}  \\
		\midrule
		1  & 0.0811  &0.8280  & 5k    & 1k(36k)  \\
		2 & 0.0109 & 0.8400    &1.2k  & 0.8 (90k)  \\
		3 & 0.0021& 0.9250&1k & 0.2 (270k) \\
		test$^1$ & 0.0004  &  &   &100k\\
		gen.$^2$ &0  &1.0000  &  & 100k\\
		\bottomrule
	\end{tabular}
	
	$^1$ denotes testing results within the training data sampling space.
	
	$^2$ denotes testing results out of the training data sampling space.
	\vspace{-5mm}
\end{table}

\begin{figure}[t]
	\centering
	\vspace{-1mm}
	\subfigure[Trajectories case 1]{
		\begin{minipage}[t]{0.22\textwidth}
			\centering
			\includegraphics[width=\textwidth]{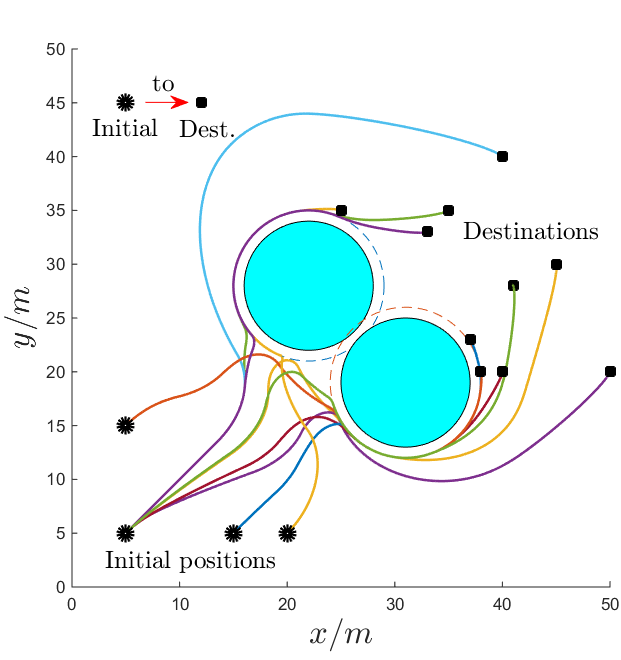}
		\end{minipage}%
	}%
	\subfigure[Trajectories case 2]{
		\begin{minipage}[t]{0.22\textwidth}
			\centering
			\includegraphics[width=\textwidth]{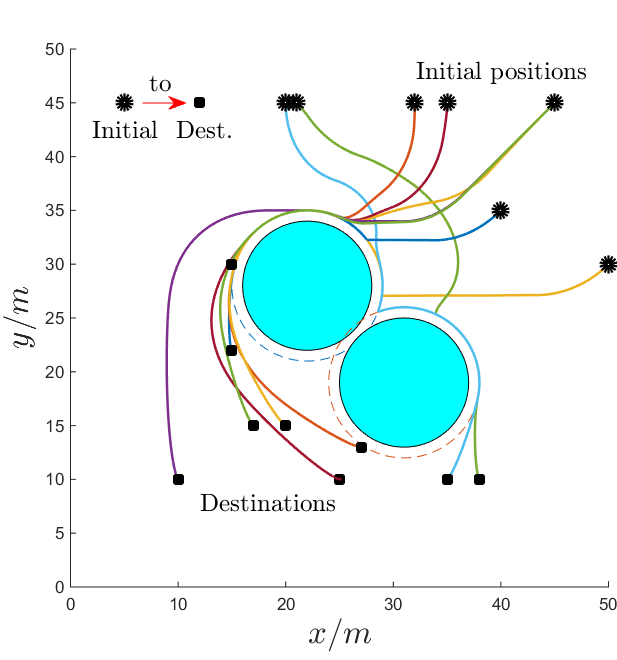}
		\end{minipage}%
	}%
	
	\centering
	\caption{Robot control problem feasibility and reachability test after learning in irregular obstacle case.}\label{fig:app}
	\vspace{-4mm}
\end{figure}

As shown in Table \ref{tab:accu2}, the classification accuracy and the infeasibility rate change similarly to the regular obstacle case in each iteration. The training results in this case are shown in Table \ref{tab:accu2}.
We also apply this hypersurface to the robot control problem, and test the feasibility and reachability. The QPs (\ref{eqn:objf}) are always feasible on the path from the initial positions to destinations. The obstacles are safely avoided, and the robot can reach the destinations in most cases. We present the results in Fig. \ref{fig:app}(a)-(b). The robot can safely avoid the local traps formed by these two circle obstacles, and thus, reachability is also improved in addition to feasibility.

\subsection{Autonomous Driving}

In autonomous driving, the ego vehicle treats all the other vehicles as moving obstacles. We consider the same dynamics and constraints as in (\ref{eqn:robot})-(\ref{eqn:robotcontrol1}), except relaxing the maximum speed limit to $28m/s$. We use the sampling-based learning approach to recursively improve the QP feasibility with respect to moving obstacles.  As all the vehicle types (such as size) are known, we can learn a feasibility constraint for each type of vehicle, and then apply this feasibility constraint to the QP. We fully cover the other vehicle by a disk, and only consider the distance between the center of the ego vehicle and the disk. 

\begin{figure}[htbp]
	\centering
	\vspace{-5mm}
	\subfigure[Snapshot time 1s]{
		\begin{minipage}[t]{0.23\textwidth}
			\centering
			\includegraphics[width=\textwidth]{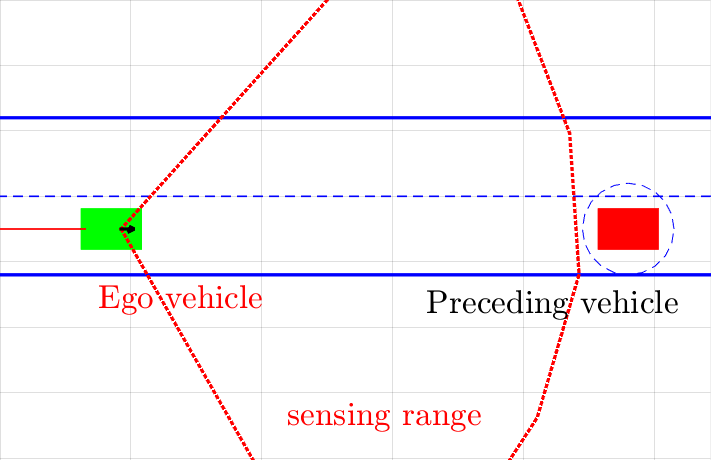}
		\end{minipage}%
	}%
	\subfigure[Snapshot time 3s]{
		\begin{minipage}[t]{0.23\textwidth}
			\centering
			\includegraphics[width=\textwidth]{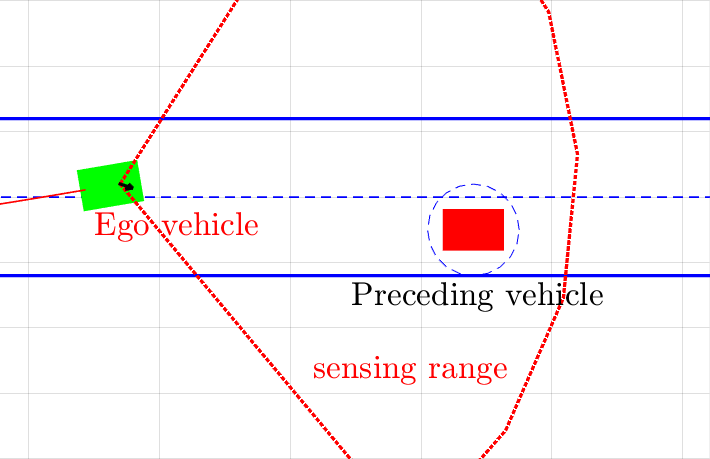}
		\end{minipage}%
	}%

    \subfigure[Snapshot time 7s]{
    	\begin{minipage}[t]{0.23\textwidth}
    		\centering
    		\includegraphics[width=\textwidth]{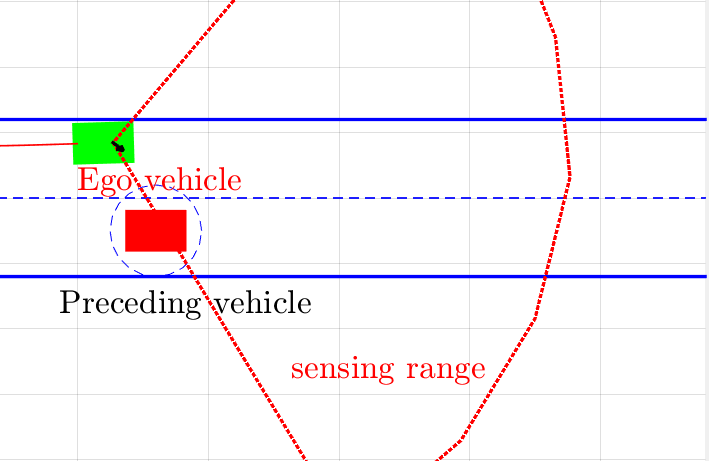}
    	\end{minipage}%
    }%
	\subfigure[Snapshot time 10s]{
		\begin{minipage}[t]{0.23\textwidth}
			\centering
			\includegraphics[width=\textwidth]{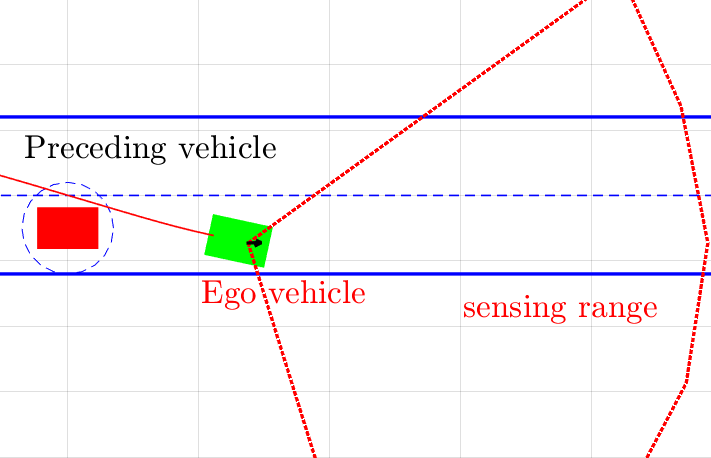}
		\end{minipage}%
	}%
	
	\centering
	\caption{Sampling learning approach applied in autonomous driving for the overtaking of another moving vehicle.}\label{fig:av}
	\vspace{-1mm}
\end{figure}

In the learning process, we take the relative position and relative speed between the ego vehicle and the other vehicle in both along-lane and lateral directions, and the heading of the ego vehicle as the inputs for the SVM model (\ref{eqn:kernel}). The relative speed difference is sampled between $0m/s$ and $20m/s$. The feedback learning process is similar to Table \ref{tab:accu}, but takes 6 iterations in order to achieve an infeasible rate that is smaller than $\varepsilon$. In the vehicle overtaking test, we set the initial and desired speeds for the ego vehicle as $28m/s$, while the other vehicle runs at a constant speed $16m/s$. Note that the control bounds for both $u_1, u_2$ are very tight, as shown in (\ref{eqn:robotcontrol1}). Without the feasibility constraint, the QP will be infeasible at some time instant for the ego vehicle. However, the ego vehicle can successfully overtake the other vehicle and the QP is always feasible with the learned feasibility constraint, as the snapshot shown in Fig. \ref{fig:av}.

\section{Conclusion}
\label{sec:conclusion}

We proposed a machine learning technique for the CBF-based QPs corresponding to optimal control problems with safety constraints and control limitations.  The learning approach can deal with both regular and irregular unsafe sets. The simulation results on a robot control problem and an autonomous driving case study show good feasibility performance with the proposed approaches. Future work will focus on dynamics and datasets affected by noise.


\bibliographystyle{IEEEtran}
\bibliography{CBF}

\end{document}